  \documentclass[11pt]{article}
  \begin{document}

\newcommand{\bm}[1]{\mbox{\boldmath $#1$}}

\def\xiv{\vec \xi }
\def\etav{\vec \eta}
\def\lie{{\pounds}}
\def\xif{\mbox{\boldmath $ \xi $}}
\def\etaf{\mbox{\boldmath $ \eta $}}
\def\t{\tau}
\def\f{\varphi}
\def\di{{\rm div}}

\title{Trapped submanifolds in Lorentzian geometry}%
\author{Jos\'e M. M. Senovilla\\
F\'{\i}sica Te\'orica, Universidad del Pa\'{\i}s Vasco\\
Apartado 644, 48080 Bilbao, Spain}%
\maketitle

\begin{abstract}
In Lorentzian geometry, the concept of trapped submanifold is
introduced by means of the mean curvature vector properties. 
Trapped submanifolds are generalizations of the standard maximal hypersurfaces 
and minimal surfaces, of geodesics, and also of the trapped surfaces introduced 
by Penrose. Selected applications to gravitational 
theories are mentioned. 
\end{abstract}


53B25,53B30,53B50

\section{Introduction}
The concept of closed trapped surface, first 
introduced by Penrose \cite{P}, is extremely useful in many 
physical problems and mathematical developments, with truly versatile applications. It was a cornerstone for the achivement of the singularity 
theorems, the analysis of gravitational collapse, 
the study of the cosmic censorship hypothesis, 
or the numerical evolution of initial data, just to mention a few, see e.g. \cite{P,HE,S1}
(a more complete list of references can be found in \cite{MS}.) Trapped
surfaces are usually introduced as co-dimension 2 imbedded spatial surfaces such that
all its local portions have, at least initially, a decreasing (increasing) area along 
{\it any} future evolution direction. However, it has been seldom recognized that the concept of trapped surface is genuinely and purely geometric, closely related to the traditional concepts of geodesics, minimal surfaces and variations of submanifolds. The purpose of this short note is to present this novel view, which may be clarifying for, and perhaps arouse interest of, the mathematical community.

\section{Basics on semi-Riemannian submanifolds.}
Let $(V,g)$ be any $D$-dimensional semi-Riemannian manifold with 
metric tensor $g$ of any signature. An imbedded submanifold is a 
pair $(S,\Phi)$ where $S$ is a $d$-dimensional manifold on its own 
and $\Phi : S \longrightarrow V$ is an imbedding \cite{O}. As is 
customary in mathematical physics, for the sake of brevity 
$S$ will be identified with its image $\Phi (S)$ in $V$. $D-d$ is 
called the co-dimension of $S$ in $V$.

At any $p\in \Phi(S)$ one has the decomposition of the tangent space
$$
T_{p}V=T_{p}S\oplus T_{p}S^{\perp}
$$
{\em if and only if} the inherited metric (or first fundamental 
form) $\Phi^{*}g\equiv \gamma$ is non-degenerate at $p$. Henceforth, I 
shall assume that $\gamma$ is non-degenerate everywhere. Let us note in 
passing that $\Phi(S)$ is called {\em spacelike} if $\gamma$ is also 
positive definite. Thus, $\forall p\in S,\,\,\, \forall \vec v\in T_{p}V$ we have $\vec v 
=\vec{v}^{T} + \vec{v}^{\perp}$ which are called the {\em tangent} 
and {\em normal} parts of $\vec v$ relative to $S$.

Obviously, $(S,\gamma)$ is a semi-Riemannian manifold on its own, and 
its intrinsic structure as such is inherited from $(V,g)$. However, 
$(S,\gamma)$ inherits also {\em extrinsic} properties. Important 
inherited intrinsic objects are (i) the canonical 
volume element $d$-form $\etaf_{S}$ associated to $\gamma$; (ii) a Levi-Civita 
connection $\overline\nabla$ such that $\overline\nabla\gamma =0$. An equivalent 
interesting characterization is
\begin{equation}
\forall \vec x,\vec y \in TS , \hspace{1cm} 
\overline\nabla_{\vec x}\, \vec y=\left(\nabla_{\vec x}\, \vec y\right)^{T} 
\label{tangent}
\end{equation}
(where $\nabla$ is the connection on $(V,g)$); and (iii) of course, 
the curvature of $\overline\nabla$ and all derived objects thereof.

Concerning the extrinsic structure, the basic object is the {\em 
shape tensor} $K : TS \times TS \longrightarrow TS^{\perp}$, also 
called extrinsic curvature of $S$ in $V$, defined by
\begin{equation}
\forall \vec x,\vec y \in TS ,\hspace{1cm}  
K(\vec x,\vec y)=-\left(\nabla_{\vec x}\, \vec y\right)^{\perp} \, .
\label{shape}
\end{equation}
The combination of (\ref{tangent}) and (\ref{shape}) provides
$$
\forall \vec x,\vec y \in TS , \hspace{1cm} 
\nabla_{\vec x}\, \vec y = \overline\nabla_{\vec x}\, \vec y
- K(\vec x,\vec y) \, .
$$
An equivalent way of expressing the same is
$$
\forall \bm{\omega}\in T^{*}S , \hspace{1cm} 
\Phi^{*}(\bm{\nabla\omega})=\overline\nabla 
(\Phi^{*}\bm{\omega})+\bm{\omega} (K)
$$
where by definition $\bm{\omega} (K) (\vec x,\vec y)=\bm{\omega} \left(K(\vec 
x,\vec y)\right)$ for all $\vec x,\vec y \in TS$.

The shape tensor contains the information concerning the ``shape'' 
of $\Phi(S)$ within $V$ along {\em all} directions normal to 
$\Phi(S)$. Observe that $K(\vec x,\vec y)\in TS^{\perp}$. If one 
chooses a particular normal direction $\vec n\in TS^{\perp}$, then one 
defines a 2-covariant symmetric tensor field $K_{\vec n}\in 
T_{(0,2)}S$ by means of
$$
K_{\vec n}(\vec x,\vec y)=\bm{n}(K)(\vec x,\vec y)=
g\left(\vec n,K(\vec x,\vec y)\right) , \hspace{3mm} 
\forall \vec x,\vec y \in TS
$$
which is called the {\em second fundamental form} of $S$ in $(V,g)$ 
relative to $\vec n$.

\section{The mean curvature vector.}
The main object to be used in this contribution is the {\em mean 
curvature vector} $\vec H$ of $S$ in $(V,g)$. This is an averaged 
version of the shape tensor defined by
$$
\vec H =\mbox{tr}\,  K , \hspace{5mm} \vec H \in TS^{\perp}
$$
where the trace tr is taken with respect to $\gamma$, of course. Each 
component of $\vec H$ along a particular normal direction, that is to 
say, $g(\vec H, \vec n)$ (= tr $K_{\vec n}$) is termed ``expansion 
along $\vec n$'' in some physical applications.

The classical interpretation of $\vec H$ can be 
understood as follows. Let us start with the simplest case $d=1$, so that 
$S$ is a curve in $V$. Then there is only one independent tangent 
vector, say $\vec x$, and $\left(\nabla_{\vec x}\, \vec 
x\right)^{\perp}=-K =- \vec H$ is simply (minus) the proper 
acceleration vector of the curve. In other words, $S$ is a geodesic if 
and only if $K=0$ (equivalently in this case, $\vec H =\vec 0$). 
Hence, an immediate and standard generalization of a geodesic to 
arbitrary codimension $d$ is: ``$S$ is {\em totally geodesic} if and only 
if $K=0$''. Totally geodesic submanifolds are those such that all 
geodesics within $(S,\gamma)$ are geodesics on $(V,g)$.

Nevertheless, one can also generalize the concept of geodesic to 
arbitrary $d$ by assuming just that $\vec H=\vec 0$. To grasp the 
meaning of this condition, let us first consider the opposite extreme 
case: $d=D-1$ or codimension 1. Then, $S$ is a hypersurface and there exists 
only one independent normal direction, say $\vec n$, so that 
necessarily $\vec H =\theta \vec n$ where $\theta$ is the (only) 
expansion, or divergence. Classical results 
imply that the vanishing of $\vec H$ (ergo $\theta =0$) defines the situation 
where there is no local variation of volume along the normal direction. 
Actually, this interpretation remains valid for arbitrary $d$. Indeed, 
let $\xiv$ be an arbitrary $C^1$ vector field on $V$ defined on a neighbourhood 
of $S$, and let $\{\f_{\t}\}_{\t\in I}$ be its flow, that is
its local one-parameter group of local transformations, where $\t$ 
is the canonical parameter and $I\ni 0$ is a real interval. This 
defines a one-parameter family of surfaces $S_{\t}\equiv \f_{\t}(S)$ 
in $V$, with corresponding imbeddings 
$\Phi_{\t}: S\rightarrow V$ given by $\Phi_{\t}=\f_{\t}\circ 
\Phi$. Observe that $S_{0}=S$. Denoting by $\etaf_{S_{\t}}$
their associated canonical volume element $d$-forms,
it is a matter of simple calculation to get
$$
\left.\frac{d\etaf_{S_{\t}}}{d\t}\right|_{\t=0}=
\frac{1}{2}\mbox{tr}\left[\Phi^{*}(\lie_{\xiv} \, g)\right]
\etaf_{S}
$$
where $\lie_{\xiv}$ is the Lie derivative with respect to $\xiv$. 
Another straightforward computation using the standard formulae 
relating the connections on $\nabla$ and $\overline\nabla$ leads to
\begin{equation}
\frac{1}{2}\mbox{tr}\left[\Phi^{*}(\lie_{\xiv} \, g)\right]=
\di (\f^{*}{\xif}) + g(\xiv , \vec{H}) \label{new}
\end{equation}
where $\di$ is the divergence operator on $S$. Combining the two 
previous formulas one readily gets the expression for the variation of 
$d$-volume:
$$
\left.\frac{dV_{S_{\t}}}{d\t}\right|_{\t=0}=\int_{S}\left(\di \vec{\bar{\xi}} + 
g(\xiv , \vec{H})\right)\, \etaf_{S} 
$$
where $V_{S_{\t}}=\int_{S_{\t}}\etaf_{S_{\t}}$ is the volume of $S_{\t}$. 
In summary: 

\vspace{2mm}
\noindent
{\em Among the set of all submanifolds without boundary 
(or with a fixed boundary under appropriate restrictions) those of 
extremal volume must have $\vec H =\vec 0$}.

\section{Lorentzian case. Future-trapped submanifolds.}
If $(V,g)$ is a proper Riemannian manifold, then $g(\vec H , \vec 
H)\geq 0$ and the only distinguished case is $g(\vec H , \vec H)= 0$ 
which is equivalent to $\vec H =\vec 0$: a extremal submanifold. 
However, in general semi-Riemannian manifolds $g(\vec H , \vec H)$ can 
be also negative, as well as zero with non-vanishing $\vec H$. Thus, 
new possibilities and distinguished cases arise.

To fix ideas, let us concentrate in the physically relevant case of 
a Lorentzian manifold $(V,g)$ with signature (--,+,\dots ,+). Let 
$(S,\gamma)$ be spacelike. Then, $\vec H$ can be classified according 
to its causal character:
$$
g(\vec H , \vec H)=\left\{
\begin{array}{cl}
>0 & \vec H\,\, \mbox{is spacelike}\\
=0 & \vec H\,\, \mbox{is null (or zero)}\\
<0 & \vec H\,\, \mbox{is timelike}
\end{array}
\right.
$$
Of course, this sign can change from point to point of $S$. 
Recall that non-spacelike vectors can be subdivided into future- 
and past-pointing. Hence, $S$
can be classified as (omitting past duals) \cite{MS,S}:
\begin{enumerate}
\item {\em future trapped} if $\vec H$ is timelike and 
future-pointing all over $S$. 
\item {\em nearly future trapped} if $\vec H$ is non-spacelike and 
future-pointing all over $S$, and timelike at least at a point of $S$.
\item {\em marginally future trapped} if $\vec H$ is null and future-pointing 
all over $S$, and non-zero at least at a point of $S$.
\item {\em extremal} or {\em symmetric} if $\vec H =\vec 0$ all over $S$.
\item {\em absolutely non-trapped} if $\vec H$ is spacelike all over $S$.
\end{enumerate}
The original definition of ``closed trapped surface'', which is of 
paramount importance in General Relativity ($D=4$), is due to 
Penrose \cite{P,HE,S1} and was for codimension two, in which case 
points 1, 4 and 5 coincide with the standard 
nomenclature; point 2 was coined in \cite{MS}, while 
3 is more general than the standard concept in GR (e.g. 
\cite{HE,S1}) ---still, all standard marginally trapped $(D-2)$-surfaces are 
included in 3---. On the other hand, the above terminology is 
unusual for the cases $d=D-1$ or $d=1$, see \cite{MS,S} for 
explanations. 

\section{Applications}
One of the advantages of having defined trapped submanifolds via $\vec 
H$ is ---apart from being generalizable to arbitrary codimension 
and thereby comparable with well-known cases such as maximal hypersurfaces 
and geodesics--- that many simple results and applications can be 
derived. As an example, let us consider the case in which $\xiv$ is a 
conformal Killing vector $\lie_{\xiv} g =2\Psi g$ (including the 
particular cases of homotheties ($\Psi =$const.) and proper Killing 
vectors ($\Psi =0$)). Then formula (\ref{new}) specializes to
$\Psi \, d = \di (\f^{*}{\xif}) + g(\xiv , \vec{H})$
so that, integrating over any {\em closed} $S$ (i.e.\, compact 
without boundary) we get
$$
\int_{S} \Psi \etaf_{S} = \frac{1}{d} \int_{S}g(\xiv , 
\vec{H})\etaf_{S} \, .
$$
Therefore, if $\Psi|_{S}$ has a sign, then $g(\xiv , \vec{H})$ must 
have the same sign, clearly restricting the possibility of $\vec 
H$ being non-spacelike. For instance, if $\xiv$ is timelike, then $\vec H$ 
(if non-spacelike) must be oppositely directed to sign$(\Psi|_{S})\xiv$; in 
particular, if $\Psi =0$, then
there cannot be closed (nearly, marginally) trapped submanifolds at all 
\cite{MS,S}. Analogously, if $\xiv$ is null on $S$ and $\Psi|_{S}=0$, 
then the only possibility for a non-spacelike $\vec H$ is that the 
mean curvature vector be null and proportional to $\xiv$.

Specific consequences of the above are, for example, \cite{MS,S,S2}
\begin{itemize}
\item that in Robertson-Walker spacetimes (where there is a conformal Killing 
vector), closed spacelike geodesics are forbiden (!), and closed 
submanifolds can only be {\em past}-trapped if the model is expanding 
\cite{MS,S}; furthermore, there cannot be maximal closed hypersurfaces, nor 
minimal surfaces \cite{MS,S}.
\item in stationary regions of $(V,g)$, any marginally trapped, 
nearly trapped, or trapped submanifold is necessarily non-closed 
and non-orthogonal to the timelike Killing vector \cite{MS,S}.
\item in regions with a null Killing vector $\xiv$, all trapped or
nearly trapped submanifolds must be non-closed and 
non-orthogonal to $\xiv$, and any marginally trapped submanifold
must have a mean curvature vector parallel (and orthogonal!)
to the null Killing vector.
\item the impossibility of existence of closed trapped surfaces 
(co-dimension 2) in spacetimes (arbitrary dimension) with vanishing curvature 
invariants \cite{S2}. This includes, in particular, the case of 
pp-waves \cite{MS,S,S2}. This has applications to modern 
string theories, implying that the spacetimes with vanishing curvature 
invariants, which are in particular exact solutions of the full 
non-linear theory, do not posses any horizons.
\end{itemize}

More details and applications can be found in \cite{MS,S,S2,G}.

\section*{Acknowledgments}
This work has been partially supported by
grant no. 9/UPV 00172.310-14456/2002 of the University of the Basque 
Country.


\begin{thebibliography}{9}
\bibitem{P} R. Penrose, Gravitational collapse and space-time singularities,
{\em Phys. Rev. Lett.} {\bf 14} (1965) 57.

\bibitem{HE} S.W. Hawking and G.F.R. Ellis, {\it The large scale
structure of space-time\/}, (Cambridge Univ. Press), Cambridge (1973).

\bibitem{S1} J.M.M. Senovilla, Singularity theorems and their 
consequences, \emph{Gen. Rel. Grav.} {\bf 30} (1998) 701.

\bibitem{MS} M. Mars and J.M.M. Senovilla, {\em Class. Quantum Grav.} 
{\bf 20} (2003) L293.

\bibitem{O} B. O'Neill, {\em Semi-Riemannian Geometry} (Academic P.), New 
York (1983)

\bibitem{S} J.M.M. Senovilla, Novel results on trapped surfaces,
in ``Mathematics of Gravitation II", 
(Warsaw, September 1-9, A Kr\'olak and K Borkowski eds, 2003);
{\tt http://www.impan.gov.pl/Gravitation/ConfProc/index.html} 
(gr-qc/03011005).

\bibitem{S2} J.M.M. Senovilla, On the existence of horizons in spacetimes 
with vanishing curvature invariants, J. High Energy Physics {\bf 11} 
(2003) 046.

\bibitem{G} D. Gerber, Trapped submanifolds in spacetimes with 
symmetries, Diploma thesis, ETH Zurich and Univ. Neuch\^atel, (2004)
{\tt http://www.unine.ch/phys/string/DiplomDG.pdf}

\end{thebibliography}
  \end{document}